\input amstex
\documentstyle{amsppt}
\magnification=\magstephalf \hsize = 6.5 truein \vsize = 9 truein
\vskip 3.5 in

\rightheadtext\nofrills {TRANSLATION EQUIVALENT ELEMENTS IN FREE GROUPS}
\NoBlackBoxes
\TagsAsMath

\def\label#1{\par%
        \hangafter 1%
        \hangindent .75 in%
        \noindent%
        \hbox to .75 in{#1\hfill}%
        \ignorespaces%
        }

\newskip\sectionskipamount
\sectionskipamount = 24pt plus 8pt minus 8pt
\def\sectionskip{\vskip\sectionskipamount}
\define\sectionbreak{%
        \par  \ifdim\lastskip<\sectionskipamount
        \removelastskip  \penalty-2000  \sectionskip  \fi}
\define\section#1{%
        \sectionbreak   
        \subheading{#1}%
        \bigskip
        }

\redefine\qed{{\unskip\nobreak\hfil\penalty50\hskip2em\vadjust{}\nobreak\hfil
    $\square$\parfillskip=0pt\finalhyphendemerits=0\par}}

        
        \let    \< = \langle
        \let    \> = \rangle

\define\op#1{\operatorname{\fam=0\tenrm{#1}}} 
        \define         \a              {\alpha}
        \redefine       \b              {\beta}
        \redefine       \d              {\delta}
        \redefine       \D              {\Delta}
        \define         \e              {\varepsilon}
        \define         \E              {\op {E}}
        \define         \g              {\gamma}
        \define         \G              {\Gamma}
        \redefine       \l              {\lambda}
        \redefine       \L              {\Lambda}
        \define         \n              {\nabla}
        \redefine       \var            {\varphi}
        \define         \s              {\sigma}
        \redefine       \Sig            {\Sigma}
        \redefine       \t              {\tau}
        \define         \th             {\theta}
        \redefine       \O              {\Omega}
        \redefine       \o              {\omega}
        \define         \z              {\zeta}
        \define         \k              {\kappa}
        \redefine       \i              {\infty}
        \define         \p              {\partial}
        \define         \vsfg           {\midspace{0.1 truein}}

\topmatter
\title Translation equivalent elements in free groups
\endtitle
\author Donghi Lee
\endauthor

\address {Department of Mathematics, Pusan National University, Jangjeon-Dong, Geumjung-Gu, Pusan 609-735, Korea}
\endaddress

\email {donghi\@pusan.ac.kr}
\endemail


\abstract {Let $F_n$ be a free group of rank $n \ge 2$. Two elements $g, h$  in $F_n$ are said to be translation equivalent in $F_n$ if the cyclic length of $\phi(g)$ equals the cyclic length of $\phi(h)$ for every automorphism $\phi$ of $F_n$. Let $F(a, b)$ be the free group generated by $\{a, b\}$ and let $w(a,b)$ be an arbitrary word in $F(a,b)$. We prove that $w(g, h)$ and $w(h, g)$ are translation equivalent in $F_n$  whenever $g, h \in F_n$ are translation equivalent in $F_n$, which hereby gives an affirmative solution to problem F38b in the online version (http://www.grouptheory.info) of [1].}
\endabstract
\endtopmatter


\document
\baselineskip=24pt

\heading 1. Introduction
\endheading

Throughout this paper, let $F_n$ be the free group of rank $n \ge 2$ on the set $\Sigma$. As usual, for a word $v$ in $F_n$, $|v|$ denotes the length of the reduced word over $\Sigma$ representing $v$. A word $v$ is called {\it cyclically reduced} if all its cyclic permutations are reduced. A {\it cyclic word} is defined to be the set of all cyclic permutations of a cyclically reduced word. By $[v]$ we denote the cyclic word associated with a word $v$. Also by $||v||$ we mean the length of the cyclic word $[v]$ associated with $v$, that is, the number of cyclic permutations of a cyclically reduced word which is conjugate to $v$. The length $||v||$ is called the {\it cyclic length} of ~$v$.

Motivated by the notion of hyperbolic equivalence of homotopy classes of closed curves on surfaces studied by Leininger [3], Kapovich-Levitt-Schupp-Shpilrain [2] introduced and studied in detail the notion of translation equivalence in free groups. The following definition is a combinatorial version of translation equivalence:

\proclaim{Definition 1.1 {\rm [2, Corollary 1.4]}} Two words $g, \ h \in F_n$ are called {\it translation equivalent in $F_n$} if the cyclic length of $\phi(g)$ equals the cyclic length of $\phi(h)$ for every automorphism $\phi$ of $F_n$.
\endproclaim

Kapovich-Levitt-Schupp-Shpilrain [2] also produced two different sources of translation equivalence in free groups. Let $F(a, b)$ be the free group with basis $\{a, b\}$.

\proclaim{Theorem 1.2 {\rm [2, Theorem C]}} Let $w(a, b) \in F(a, b)$ be a freely reduced word. Then for any $g, \, h \in F_n$, $w(g, h)$ and $w^R(g, h)$ are translation equivalent in $F_n$, where $w^R(a, b)=w(a^{-1}, b^{-1})^{-1}$.
\endproclaim

\proclaim{Theorem 1.3 {\rm [2, Theorem D]}} Let $g, \, h \in F_n$ be translation equivalent in $F_n$ but $g \neq h^{-1}$. Then for any positive integers $p, q, i, j$ such that $p+q=i+j$, $g^ph^q$ and $g^ih^j$ are translation equivalent in $F_n$.
\endproclaim

Providing these sources, Kapovich-Levitt-Schupp-Shpilrain [2] asked whether or not $w(g, h)$ and $w(h, g)$ are translation equivalent in $F_n$ whenever $g, \, h \in F_n$ are translation equivalent in $F_n$ and $w(a, b) \in F(a, b)$ is arbitrary. The purpose of the present paper is to answer this question in the affirmative:

\proclaim{Theorem 1.4} Let $g, \, h \in F_n$ be translation equivalent in $F_n$ and let $w(a, b) \in F(a, b)$ be arbitrary. Then $w(g, h)$ and $w(h, g)$ are translation equivalent in $F_n$.
\endproclaim

The proof of this theorem is quite combinatorial and will appear in the next section. We remark that the same technique used in the present paper can also be used to prove Theorem~1.3, thus giving an alternative proof of the theorem. In [2] Theorem~1.3 was proved by a geometric method that uses the analysis of possible axis configurations for compositions of isometries of $\Bbb{R}$-trees.

\heading 2. Proof of Theorem 1.4
\endheading

If the cyclic word $[w(a, b)]$ associated with $w(a, b)$ is trivial, then there is nothing to prove. So let the cyclic word $[w(a, b)]$ be nontrivial. If $g= h^{\pm 1}$, then the statement clearly holds. So suppose that $g \neq h^{\pm 1}$. If the cyclic word $[w(a, b)]$ consists of only one letter, that is, $[w(a, b)]=[a^i]$ or $[w(a, b)]=[b^i]$ for some nonzero $i$, then the statement is obvious, because $g^i$ and $h^i$ are translation equivalent in $F_n$. Otherwise, the cyclic word $[w(a, b)]$ can be written as
$$[w(a, b)]=[a^{\ell_1} b^{m_1} \cdots a^{\ell_k} b^{m_k}],
\tag 2.1
$$
where $\ell_i, \, m_i$ are nonzero integers. Let $\phi$ be an arbitrary automorphism of $F_n$. To prove the theorem, we need to show that $||\phi(w(g, h))||=||\phi(w(h, g))||$. For simplicity, we write
$$X=\phi(g) \quad \text{and} \quad Y=\phi(h).$$
The words $X$ and $Y$ in $F_n$ can be uniquely factored as
$$X=A \bar{X} A^{-1} \ \text{(reduced)} \quad \text{and} \quad Y=B \bar{Y} B^{-1}\ \text{(reduced)},
\tag 2.2
$$
where $\bar{X}$ and $\bar{Y}$ are cyclically reduced ($A$ or $B$ may be the empty word). It then follows from the translation equivalence of $g$ and $h$ in $F_n$ that $$|\bar{X}|=|\bar{Y}|.
\tag 2.3
$$
Combining (2.1) with (2.2) gives
$$\aligned
[\phi(w(g, h))]&=[A \bar{X}^{\ell_1} A^{-1} B \bar{Y}^{m_1} B^{-1}  \cdots A \bar{X}^{\ell_k} A^{-1} B \bar{Y}^{m_k} B^{-1}] \\
[\phi(w(h, g))]&=[B \bar{Y}^{\ell_1} B^{-1} A \bar{X}^{m_1} A^{-1}  \cdots B \bar{Y}^{\ell_k} B^{-1} A \bar{X}^{m_k} A^{-1}].
\endaligned
\tag 2.4
$$

We continue to argue dividing into three cases according to which part is cancelled in the product $A^{-1}B$.

\proclaim {Case 1} Neither $A^{-1}$ nor $B$ is completely cancelled in the product $A^{-1}B$.
\endproclaim

Let $\ell=\sum_{i=1}^k |\ell_i|$ and $m=\sum_{i=1}^k |m_i|$. In this case (2.4) yields that
$$\aligned
||\phi(w(g, h))||&=\ell |\bar{X}|+m|\bar{Y}|+k|A^{-1}B|+k|B^{-1}A| \\
||\phi(w(h, g))||&=\ell |\bar{Y}|+m|\bar{X}|+k|B^{-1}A|+k|A^{-1}B|,
\endaligned$$
so by (2.3) $||\phi(w(g, h))||=||\phi(w(h, g))||$, as desired.

\proclaim {Case 2} Both $A^{-1}$ and $B$ are completely cancelled in the product $A^{-1}B$, that is, $A=B$.
\endproclaim

In this case, we have from (2.4) that
$$\aligned
[\phi(w(g, h))]&=[\bar{X}^{\ell_1} \bar{Y}^{m_1} \cdots \bar{X}^{\ell_k} \bar{Y}^{m_k}] \\
[\phi(w(h, g))]&=[\bar{Y}^{\ell_1} \bar{X}^{m_1} \cdots \bar{Y}^{\ell_k} \bar{X}^{m_k}].
\endaligned
\tag 2.5
$$
At this point we prove the following

\proclaim {Claim} In the product of the form $\bar{X}^{\epsilon_1} \bar{Y}^{\epsilon_2} \bar{X}^{\epsilon_3}$, the middle part $\bar{Y}^{\epsilon_2}$ cannot be completely cancelled for any $\epsilon_i=\pm 1$.
\endproclaim

\demo{Proof of the claim} Suppose to the contrary that $\bar{Y}^{\epsilon_2}$ is completely cancelled in the product $\bar{X}^{\epsilon_1} \bar{Y}^{\epsilon_2} \bar{X}^{\epsilon_3}$. If $\bar{X}^{\epsilon_1}= \bar{Y}^{-\epsilon_2}$ or $\bar{X}^{\epsilon_3}=\bar{Y}^{-\epsilon_2}$, then we have $\bar{X}=\bar{Y}^{\pm 1}$. It then follows from $A=B$ that
$X=A\bar{X}A^{-1}=B\bar{Y}^{\pm 1}B^{-1}=Y^{\pm 1}$, a contradiction to our assumption that $g \neq h^{\pm 1}$. Otherwise, the only possibility is by (2.3) that $\bar{X}^{\epsilon_1}$ ends with  $\bar{Y_1}^{-\epsilon_2}$, $\bar{X}^{\epsilon_3}$ begins with $\bar{Y_2}^{-\epsilon_2}$, and that $\bar{Y}^{\epsilon_2}=\bar{Y_1}^{\epsilon_2} \bar{Y_2}^{\epsilon_2}$, where $|\bar{Y_i}|>1$ for all $i=1, \, 2$. Since $\bar{Y}$ is cyclically reduced, this can happen only when $\epsilon_1=\epsilon_3$. This together with (2.3) yields that $\bar{X}^{\epsilon_1}=\bar{Y_2}^{-\epsilon_2} \bar{Y_1}^{-\epsilon_2}=\bar{Y}^{-\epsilon_2}$. But then again we have $\bar{X}=\bar{Y}^{\pm 1}$. This contradiction completes the proof of the claim.
\qed
\enddemo

Similarly it can be shown that in the product of the form $\bar{Y}^{\epsilon_1} \bar{X}^{\epsilon_2} \bar{Y}^{\epsilon_3}$, the middle part $\bar{X}^{\epsilon_2}$ cannot be completely cancelled for any $\epsilon_i=\pm 1$.

For $x, \, y \in \{a, b\}^{\pm 1}$, let $n(xy)$ denote the number of occurrences of the subword $xy$ in the cyclic word $[w(a, b)]$. We put
$$\aligned
\alpha_1=n(ab) \quad &\text{and} \quad \alpha_2=n(ba); \\
\beta_1=n(a^{-1}b^{-1}) \quad &\text{and} \quad \beta_2=n(b^{-1}a^{-1}); \\
\gamma_1=n(ab^{-1}) \quad &\text{and} \quad \gamma_2=n(ba^{-1});\\
\delta_1=n(a^{-1}b) \quad &\text{and} \quad \delta_2=n(b^{-1}a).
\endaligned
$$
Then it is not hard to see that
$$\aligned
\alpha_1+\beta_1+\gamma_1+\delta_1&=\alpha_2+\beta_2+\gamma_2+\delta_2, \\
2\beta_1+\gamma_1+\delta_1&=2\beta_2+\gamma_2+\delta_2.
\endaligned
\tag 2.6
$$
Put $\kappa=(\gamma_1+\delta_1)-(\gamma_2+\delta_2)$ ($\kappa$ may be zero). The second equality in (2.6) gives $\kappa=(\gamma_1+\delta_1)-(\gamma_2+\delta_2)=2\beta_2-2\beta_1$. Hence $\kappa$ is an even integer and  $\beta_2=\beta_1 + \kappa/2$. This together with the first equality in (2.2) yields $\alpha_2=\alpha_1 + \kappa/2$, so that $$\alpha_1+\beta_2=\alpha_2+\beta_1.
\tag 2.7
$$

Put
$$\aligned
p_1=|\bar{X}|+|\bar{Y}|-|\bar{X}\bar{Y}| \quad &\text{and} \quad p_2=|\bar{Y}|+|\bar{X}|-|\bar{Y}\bar{X}|\\
q_1=|\bar{X}^{-1}|+|\bar{Y}^{-1}|-|\bar{X}^{-1}\bar{Y}^{-1}| \quad &\text{and} \quad q_2=|\bar{Y}^{-1}|+|\bar{X}^{-1}|-|\bar{Y}^{-1}\bar{X}^{-1}|\\
r_1=|\bar{X}|+|\bar{Y}^{-1}|-|\bar{X}\bar{Y}^{-1}| \quad &\text{and} \quad r_2=|\bar{Y}|+|\bar{X}^{-1}|-|\bar{Y}\bar{X}^{-1}| \\
s_1=|\bar{X}^{-1}|+|\bar{Y}|-|\bar{X}^{-1}\bar{Y}| \quad &\text{and} \quad s_2=|\bar{Y}^{-1}|+|\bar{X}|-|\bar{Y}^{-1}\bar{X}|\\
\endaligned
$$
Then in view of (2.5) and the Claim we have
$$\aligned
||\phi(w(g, h))||&=\ell|\bar{X}|+m|\bar{Y}|-\sum_{i=1}^2 p_i\alpha_i - \sum_{i=1}^2 q_i\beta_i-\sum_{i=1}^2 r_i\gamma_i-\sum_{i=1}^2 s_i\delta_i \\
||\phi(w(h, g))||&=\ell|\bar{Y}|+m|\bar{X}|-\sum_{i=1}^2 p_j\alpha_i - \sum_{i=1}^2 q_j\beta_i-\sum_{i=1}^2 r_j\gamma_i-\sum_{i=1}^2 s_j\delta_i,
\endaligned
\tag 2.8
$$
where $j=2$ provided $i=1$; $j=1$ provided $i=2$. Note that $p_1=q_2$, $q_1=p_2$, $r_1=r_2$ and $s_1=s_2$. This together with (2.7) yields
$$\split
p_1\alpha_1+q_2\beta_2&=p_1\alpha_1+p_1\beta_2=p_1(\alpha_1+\beta_2) \\
&=q_2(\alpha_2+\beta_1)=q_2\alpha_2+q_2\beta_1=p_1\alpha_2+q_2\beta_1,
\endsplit
$$
$$\split
p_2\alpha_2+q_1\beta_1&=p_2\alpha_2+p_2\beta_1=p_2(\alpha_2+\beta_1) \\
&=q_1(\alpha_1+\beta_2)=q_1\alpha_1+q_1\beta_2=p_2\alpha_1+q_1\beta_2,
\endsplit
$$
and that
$$r_1 \gamma_1=r_2 \gamma_1, \quad \ r_2 \gamma_2=r_1 \gamma_2, \quad s_1 \delta_1=s_2 \delta_1, \quad s_2 \delta_2=s_1 \delta_2.$$
Combining all these equalities and (2.3) with (2.8), we finally have $||\phi(w(g, h))||=||\phi(w(h, g))||$,
as required.

\proclaim {Case 3} Only one of $A^{-1}$ and $B$ is completely cancelled in the product $A^{-1}B$.
\endproclaim

Without loss of generality, assume that only $A^{-1}$ is completely cancelled in the product $A^{-1}B$. Then $B$ can be factored as
$$B=AC \ (\text{reduced}),
\tag 2.9
$$ where $|C|>0$, and hence (2.4) is rewritten as
$$\aligned
[\phi(w(g, h))]&=[\bar{X}^{\ell_1} C \bar{Y}^{m_1} C^{-1}  \cdots \bar{X}^{\ell_k} C \bar{Y}^{m_k} C^{-1}], \\
[\phi(w(h, g))]&=[\bar{Y}^{\ell_1} C^{-1} \bar{X}^{m_1} C \cdots \bar{Y}^{\ell_k} C^{-1} \bar{X}^{m_k} C].
\endaligned
\tag 2.10
$$
If the right-hand expressions of (2.10) are reduced as written, the desired equality $||\phi(w(g, h))||=||\phi(w(h, g))||$ is obvious by (2.3). Otherwise, since $C \bar{Y}^i C^{-1}$ is reduced for all nonzero $i$ (note that $B \bar{Y}B^{-1}$ is reduced in (2.2)) and $\bar{X}$ is cyclically reduced, initial cancellation can occur only between $C^{-1}$ and $\bar{X}^j$ or between $\bar{X}^j$ and $C$ but not both. Suppose that initial cancellation occurs between $\bar{X}^j$ and $C$ with $j>0$ (the other case is similar). Write
$$C=DE \ (\text{reduced}),
\tag 2.11
$$
where $E$ is the part of $C$ that remains uncancelled in the product $C^{-1}\bar{X} C$ ($E$ may be the empty word). Define $\bar{Z}$ to be the reduced form of $D^{-1} \bar{X} D$. Then notice by the choice of $D$ that $\bar{Z}$ is a cyclic permutation of $\bar{X}$, so that $\bar{Z}$ is cyclically reduced and $|\bar{Z}|=|\bar{X}|$.

From (2.10) and (2.11) we have
$$\aligned
[\phi(w(g, h))]&=[\bar{Z}^{\ell_1} E \bar{Y}^{m_1} E^{-1} \cdots  \bar{Z}^{\ell_k} E \bar{Y}^{m_k} E^{-1}], \\
[\phi(w(h, g))]&=[\bar{Y}^{\ell_1} E^{-1} \bar{Z}^{m_1} E  \cdots \bar{Y}^{\ell_k} E^{-1} \bar{Z}^{m_k} E].
\endaligned
\tag 2.12
$$
First let $E$ be a non-empty word. By the choice of $E$, $E^{-1}  \bar{Z}^i E$ is reduced for all nonzero $i$. Also since $B \bar{Y}^j B^{-1}$ is reduced for all nonzero $j$ by (2.2), $E \bar{Y}^j E^{-1}$ is reduced for all nonzero $j$. These imply that the right-hand expressions of (2.12) are reduced, so that the equality $||\phi(w(g, h))||=||\phi(w(h, g))||$ follows immediately from the fact that $|\bar{Z}|=|\bar{Y}|$.
Now let $E$ be the empty word. Then in view of (2.2), (2.9) and (2.11) we get
$$X=B \bar{Z} B^{-1}.
\tag 2.13
$$
In the product of the form $\bar{Z}^{\epsilon_1} \bar{Y}^{\epsilon_2} \bar{Z}^{\epsilon_3}$, the middle part $\bar{Y}^{\epsilon_2}$ cannot be completely cancelled for any $\epsilon_i=\pm 1$, for otherwise reasoning in the same way as in the Claim of Case ~2 we would have $\bar{Z}=\bar{Y}^{\pm 1}$, that is, $X=B \bar{Z} B^{-1}= B \bar{Y}^{\pm 1} B^{-1}=Y^{\pm 1}$ by (2.13), so that $g=h^{\pm 1}$, a contradiction. Similarly the middle part $\bar{Z}^{\epsilon_2}$ in the product of the form $\bar{Y}^{\epsilon_1} \bar{Z}^{\epsilon_2} \bar{Y}^{\epsilon_3}$ cannot be completely cancelled for any $\epsilon_i=\pm 1$. Hence we can follow the proof of Case ~2 to obtain the assertion.

The proof of the theorem is now completed.
\qed

\heading Acknowledgement
\endheading

The author is grateful to the referee for a careful report. The author was supported by Pusan National University Research Grant.

\heading References
\endheading

\roster

\item"1." G. Baumslag, A. G. Myasnikov and V. Shpilrain, {\it Open
problems in combinatorial group theory}, Contemp. Math. {\bf 250}
(1999), 1--27.

\item"2." I. Kapovich, G. Levitt, P. E. Schupp and V. Shpilrain, Translation equivalence in free groups, {\it Trans. Amer. Math. Soc.}, to appear.

\item"3." C. J. Leininger, Equivalent curves in surfaces, {\it Geom. Dedicata} {\bf 102} (2003), 151--177.
\endroster
\enddocument